\input ./style/arxiv-general.cfg
\documentclass[MSNbibl,number,citesort,dvips]{arxbj}
\makeatletter
   \@ifpackageloaded{graphicx}{}{\usepackage{graphicx}}
\makeatother
\usepackage{mathbh}
\usepackage{mathrsfs}

\aid{0}
\volume{21}
\issue{4}
\pubyear{2015}
\firstpage{2139}
\lastpage{2156}
\doi{10.3150/14-BEJ638} % kopijuoti is 'New paper accepted'
\docsubty{FLA}

\makeatletter

\newcommand{\lleft}{\left}
\newcommand{\rrvert}{\vert}
\newcommand{\rright}{\right}
\newcommand{\llvert}{\vert}
\newtheorem{theorem}{Theorem}[section]
\newtheorem{lemma}[theorem]{Lemma}
\newproclaim{definition}[theorem]{Definition}
\newtheorem{corollary}[theorem]{Corollary}
\newremark{remark}[theorem]{Remark}

\newcommand{\R}{\mathbb{R}}
\newcommand{\E}{\mathbb{E}}
\newcommand{\cvlaw}{\stackrel{\mathrm{law}}{\rightarrow}}
\newcommand{\cont}[1]{\stackrel{#1}{\frown}}
\newcommand{\converge}{\displaystyle \mathop{\longrightarrow}_{n \rightarrow+\infty}}
\newcommand{\eqref}[1]{(\ref{#1})}
\newcommand{\mathds}[1]{\mathbh{#1}}

\renewcommand{\emptyset}{\varnothing}

\newcommand{\NC}{\operatorname{NC}}

\def\vfrac#1#2{(#1)/#2}

\makeatother

\begin{document}
\begin{frontmatter}

\title{Poisson convergence on the free Poisson algebra}
\runtitle{Poisson convergence on the free Poisson algebra}

\begin{aug}
%%%% inicialai - be tarpu
\author[A]{\inits{S.}\fnms{Solesne}~\snm{Bourguin}\corref{}\ead[label=e1]{solesne.bourguin@gmail.com}}% \and
%\author{\inits{}\fnms{}~\snm{}\thanksref{}\ead[label=e2]{}}
%\author{\inits{}\fnms{}~\snm{}}
%%\runauthor{} %% auto
%\dedicated{}
\address[A]{Department of Mathematical Sciences, Carnegie Mellon
University, Pittsburgh, PA 15213, USA. \printead{e1}}
%\address[]{}
\end{aug}

% HISTORY:
\received{\smonth{12} \syear{2013}}
\revised{\smonth{4} \syear{2014}}

% ABSTRACT
%
\begin{abstract}
Based on recent findings by Bourguin and Peccati, we give a \textit
{fourth moment} type condition for an element of a free Poisson chaos
of arbitrary order to converge to a free (centered) Poisson
distribution. We also show that free Poisson chaos of order strictly
greater than one do not contain any non-zero free Poisson random
variables. We are also able to give a sufficient and necessary
condition for an element of the first free Poisson chaos to have a free
Poisson distribution. Finally, depending on the parity of the
considered free Poisson chaos, we provide a general counterexample to
the naive universality of the semicircular Wigner chaos established by
Deya and Nourdin as well as a transfer principle between the Wigner and
the free Poisson chaos.
\end{abstract}

% KEYWORDS
% visi is mazosios raides ir pagal abecele
%
\end{frontmatter}

%s1 #&#
\section{Introduction and background}

%s1.1 #&#
\subsection{Overview}
Let $ \lbrace W(t) \dvtx  t \geq0  \rbrace$ be a standard
Brownian motion on $\R_{+}$ and let $q \geq1$ be an integer. Denote by
$I_{q}^{W} ( f )$ the multiple stochastic Wiener--It\^o
integral of order $q$ of a deterministic symmetric function $f \in
L^2 ( \R_{+}^{q} )$. Denote by $L^2_s ( \R_{+}^{q} )$
the subset of $L^2 ( \R_{+}^{q} )$ composed of symmetric
functions. The collection of random variables $ \lbrace
I_{q}^{W} ( f ) \dvtx  f \in L^2_s ( \R_{+}^{q} )
 \rbrace$ is what is usually called\vspace*{2pt} the $q$th Wiener chaos
associated with $W$. In a seminal paper of 2005, Nualart and Peccati
\citep{nuapec1} proved that convergence to the standard normal
distribution of an element with variance one living inside a fixed
Wiener chaos was equivalent to the convergence of the fourth moment of
this element to three. This result is now known as the \textit{fourth
moment theorem} and can be stated as follows.

%
%th1.1 #&#
\begin{theorem}[(Nualart and Peccati \cite{nuapec1})]
\label{fourthmomenttheoremwienergaussian}
Fix an integer $q \geq2$ and let $ \lbrace f_n \dvtx n \geq1
 \rbrace$ be a sequence of functions in $L^2_s ( \R
_{+}^{q} )$ such that, for each $n \geq1$, $\E (I_{q}^{W}
( f_n )  ) = q! \Vert f_n \Vert_{L^2 ( \R_{+}^{q}
)}^{2} = 1$. Then, the following two assertions are equivalent, as $n
\rightarrow\infty$:
\begin{enumerate}[(ii)]
\item[(i)] $I_{q}^{W} ( f_n ) \cvlaw\mathcal{N}(0,1)$;
\item[(ii)] $\E ( I_{q}^{W} ( f_n )^4 ) \rightarrow\E
 (\mathcal{N}(0,1)^4  )= 3$.
\end{enumerate}
\end{theorem}

This result has led to a wide collection of new results and inspired
several new research directions -- see the book \citep{np-book}, as
well as the constantly updated webpage
\begin{center}
\href{https://sites.google.com/site/malliavinstein/home}{https://sites.google.com/site/malliavinstein/home}.
\end{center}
In \citep{noupecnoncentral}, the authors obtained a similar criterion
for non-central convergence to the Gamma distribution on a Wiener
chaos. In \citep{noupecstein}, the Malliavin calculus of variations was
combined with Stein's method to obtain quantitative versions of these
fourth moment theorems. In the framework of non-commutative
probability, Kemp \textit{et al.} \citep{knps} obtained an analog of
Theorem~\ref{fourthmomenttheoremwienergaussian} for multiple integrals
with respect to a free Brownian motion.

%
%re1.2 #&#
\begin{remark}
\label{remarlmirror}
Observe that Theorem~\ref{fourthmomenttheoremwienergaussian} is stated
for symmetric functions in $L^2(\R_{+}^q)$. In free probability theory,
the symmetry assumption can be weakened to \textit{mirror symmetry} for
defining and working with free multiple integrals (see Section~\ref{subsecfreepoissonalg} and definitions therein). We say that an element
$f$ of $L^2(\R_{+}^q)$ (the collection of all complex-valued functions
on $\R_{+}^q$ that are square-integrable with respect to the Lebesgue
measure) is \textit{mirror symmetric} if $f(t_1,\ldots,t_q) = \overline
{f(t_q,\ldots,t_2,t_1)}$, for almost every vector $(t_1,\ldots,t_q)\in
\R_{+}^q$.
\end{remark}

The aim of this paper is to investigate the convergence of sequences of
multiple integrals with respect to a free Poisson measure. More
precisely, denote by $(\mathscr{A},\varphi)$ a free probability space
and let $ \lbrace\hat{N}(B) \dvtx  B \in\mathscr{B} ( \R
_{+} )  \rbrace$, where $\mathscr{B} ( \R_{+} )$
denotes the Borels sets of $\R_{+}$, be a centered free Poisson measure
on this space. For an integer $q \geq1$, the free Poisson multiple
integral of order $q$ of a mirror-symmetric bounded function with
bounded support $f \in L^2 ( \R_{+}^{q} )$ is denoted
$I_{q}^{\hat{N}} ( f )$. Random variables of this type compose
the so-called free Poisson chaos of order $q$ associated with $\hat
{N}$. The above mentioned objects are defined and constructed in
Section~\ref{section2}; refer to that section for more details. A first
result in this direction has recently been obtained by Bourguin and
Peccati in \citep{bope1}, who proved that a fourth moment type theorem
(for semicircular limits) holds on the free Poisson algebra.

%th1.3 #&#
\begin{theorem}[(Bourguin and Peccati \cite{bope1})]
\label{fourthmomenttheoremfreepoisemi}
Fix an integer $q \geq1$ and let $\{f_n \dvtx  n\geq1\} \subset L^2(\R
_+^q)$ be a { tamed} sequence of mirror symmetric kernels such that $\|
f_n\|_{L^2(\R_+^q)}^2 \to1$,
as $n\to\infty$.\vspace*{-2pt} Then, $I^{\hat{N}}_q(f_n)$ converges in law to the
semicircular distribution $\mathcal{S}(0,1)$ if and only if $\varphi
(I^{\hat{N}}_q(f_n)^4) \to2$.
\end{theorem}

The notion of \textit{tamed sequence} of kernels is introduced in
Section~\ref{section2} below: this additional assumption has been
introduced in \citep{bope1} in order to deal with the complicated
combinatorial structures arising from the computation of moments (more
precisely, it is a sufficient condition in order to preserve spectral
bounds when converging on the free Poisson algebra). Bourguin and
Peccati also proved that a \textit{transfer principle} (as the one
established in \citep{knps}, Theorem~1.8) between the classical and the
free Poisson chaos cannot hold in full generality by providing an
example where, for the same sequence of kernels, the free Poisson
multiple integral converges towards a semicircular distribution when
the corresponding classical Poisson multiple integral converges towards
a classical Poisson limit, instead of a Gaussian distribution as would
be expected in this type of transfer results.

A natural question that arises in this context is to find moment
conditions (and potentially a unique condition in the form of a linear
combination of moments as in other fourth moment type results) to
ensure the convergence of free Poisson multiple integrals to a free
Poisson distribution. One could be tempted to consider such a question
as an analog, in the free case, of the results obtained in \citep{pec1}
or \citep{bope1} for the classical case (Poisson approximations on the
classical Poisson space), but the classical analog of the free Poisson
distribution is actually the Gamma distribution as pointed out in \citep
{nope1}, Section~1.1 (see also \citep{NicSpe}, page 203). A recent
partial result in the classical framework addressing the convergence of
Poisson multiple integrals to the Gamma distribution has been obtained
in \citep{pectha1}, Theorem~2.6.

In the free case, the closest result in this direction is the fourth
moment type characterization obtained by Nourdin and Peccati in \citep
{nope1}, Theorem~1.4, for the convergence of multiple Wigner integrals
to a free Poisson distribution, this result being itself an analogue of
the main result in \citep{noupecnoncentral}, proved by the same
authors. This first free Poisson approximation theorem reads as follows.

%th1.4 #&#
\begin{theorem}[(Nourdin and Peccati \cite{nope1})]
\label{fourthmomenttheoremwignerfreepoi}
Let $f \in L^2 ( \R_{+}^{q} )$. Let $I_{q}^{S}(f)$ denote the
multiple Wigner integral of order $q$ of $f$. Let $Z(\lambda)$ have a
free centered Poisson distribution with rate $\lambda> 0$, fix an even
integer $q \geq2$ and let $ \lbrace f_n \dvtx  n \geq1
\rbrace$ be a sequence of mirror-symmetric functions in $L^2 ( \R
_{+}^{q} )$ such that, for each $n \geq1$, $\varphi
(I_{q}^{S} ( f_n )^2  ) = \Vert f_n \Vert_{L^2 ( \R
_{+}^{q} )}^{2} = \lambda$. Then, the following two assertions are
equivalent, as $n \rightarrow\infty$:
\begin{enumerate}[(ii)]
\item[(i)] $I_{q}^{S} ( f_n ) \cvlaw Z(\lambda)$;
\item[(ii)] $\varphi ( I_{q}^{S} ( f_n )^4 ) -
2\varphi ( I_{q}^{S} ( f_n )^3 ) \rightarrow\varphi
 (Z(\lambda)^4  ) - 2\varphi (Z(\lambda)^3  ) =
2\lambda^2 - \lambda$.
\end{enumerate}
\end{theorem}

The following theorem, and main result of this paper, is a fourth
moment type characterization of the convergence to free Poisson limits
on the free Poisson algebra, and could be seen as a free counterpart of
\citep{pectha1}, Theorem~2.6.

%th1.5 #&#
\begin{theorem}
\label{fourthmomenttheoremfreepoifreepoi}
Let $Z(\lambda)$ have a free centered Poisson distribution with rate
$\lambda> 0$. Fix an integer $q \geq1$ and let $ \lbrace f_n
\dvtx  n \geq1  \rbrace$ be a tamed sequence of
mirror-symmetric kernels in $L^2 ( \R_{+}^{q} )$ such that,
for each $n \geq1$, $\varphi (I_{q}^{\hat{N}} ( f_n )^2
 ) = \Vert f_n \Vert_{L^2 ( \R_{+}^{q} )}^{2} = \lambda$.
Then, the following two assertions are equivalent, as $n \rightarrow
\infty$:
\begin{enumerate}[(ii)]
\item[(i)] $I_{q}^{\hat{N}} ( f_n ) \cvlaw Z(\lambda)$;
\item[(ii)] $\varphi ( I_{q}^{\hat{N}} ( f_n )^4 ) -
2\varphi ( I_{q}^{\hat{N}} ( f_n )^3 ) \rightarrow
\varphi (Z(\lambda)^4  ) - 2\varphi (Z(\lambda)^3  )
= 2\lambda^2 - \lambda$.
\end{enumerate}
\end{theorem}

In \citep{nope1}, Nourdin and Peccati mention that one cannot expect to
have convergence of a Wigner multiple integral of odd order to the free
Poisson distribution since these integrals have all odd moments equal
to zero, as opposed to the free Poisson distribution. It is worth
pointing out that, unlike in Theorem~\ref
{fourthmomenttheoremwignerfreepoi}, even as well as odd orders can be
considered in Theorem~\ref{fourthmomenttheoremfreepoifreepoi}. The
reason for this difference will become clear in the sequel.
Additionally, Theorem~\ref{fourthmomenttheoremfreepoifreepoi} is given
in a way more general setting than the classical case result proved in
\citep{pectha1}, Theorem~2.6. Indeed, \citep{pectha1}, Theorem~2.6, was
only proved for multiple integrals of even orders, although convergence
of multiple integrals of odd orders is not excluded by the authors who
claim that the odd case is more intricate to analyze. Also, the fourth
moment characterization in \citep{pectha1}, Theorem~2.6, was only
obtained for sequences of multiple integrals of order two. The main
result in this paper, namely Theorem~\ref
{fourthmomenttheoremfreepoifreepoi}, is given for multiple integrals of
any order and the validity of the fourth moment characterization is not
restricted to the order two.

As a consequence of Theorem~\ref{fourthmomenttheoremfreepoifreepoi}, we
prove a corollary giving insights on the (almost completely unknown)
structure of the free Poisson algebra.

%co1.6 #&#
\begin{corollary}
\label{elcorollary}
The two following statements hold:
\begin{enumerate}[(2)]
\item[(1)] Let $q \geq2$ be an integer, and let $F \neq0$ be in the
$q$th free Poisson chaos. Then, $F$ cannot have a free Poisson distribution.
\item[(2)] Let $f \in L^2 (\R_{+}  )$ be such that $\Vert f
\Vert_{L^2 (\R_{+}  )}^{2} = \lambda> 0$. Then, $I_{1}^{\hat
{N}}(f)$ has a free Poisson distribution with parameter $\lambda$ if
and only if $f$ takes values in $ \lbrace0,1  \rbrace$.
\end{enumerate}
\end{corollary}

%re1.7 #&#
\begin{remark}
Point (ii) in Corollary~\ref{elcorollary} shows that the first free
Poisson chaos does not contain only free Poisson distributions, as
opposed to the first Wigner chaos that only consists of semicircular
distributed elements. The following example illustrates the situation
where an element of the first free Poisson chaos doesn't have a free
Poisson distribution. Let $A_1$ and $A_2$ be two orthogonal (with
respect to $\mu$) Borel sets of $\R_+$ such that $\mu(A_1) = \mu(A_2) =
1$. Let $Z(8)$ be a free Poisson random variable with parameter
$\lambda= 8$. Finally, let $f$ be the function of $L^2 ( \R
_{+} ) $ defined by $f = 2 (\mathds{1}_{A_1} + \mathds{1}_{A_2})$.
It is easily checked that $\varphi ( I_{1}^{\hat{N}}(f)^2  )
= 8 = \lambda$. According to the upcoming Proposition~\ref{prop2}, it
holds that
\begin{eqnarray*}
\varphi \bigl( I_{1}^{\hat{N}}(f)^3 \bigr) = \bigl(
f \star_{1}^{0} f \bigr) \cont{1} f = \int
_{\R_{+}}f^3 \,\mathrm{d}\mu= 16.
\end{eqnarray*}
As $\varphi ( Z(8)^3  ) = 8$, $I_{1}^{\hat{N}}(f)$ cannot have
a free Poisson distribution.
\end{remark}

Finally, we provide a general counterexample to the naive universality
of the free Wigner chaos (as stated in \citep{noudey1}), where it is
proven that for multilinear homogeneous sums of free random
variables, a universality phenomenon happens in the sense that it is
sufficient that a multiple Wigner integral with an appropriate
``homogeneous'' kernel converges to the semicircular distribution for
multiple integrals with respect to any free random measure (with the
same kernel) to converge as well to the semicircular distribution. The
following theorem could also be seen as a transfer principle between
the Wigner and free Poisson chaos for multiple integrals of even
orders. This theorem should be compared with the counterexample to the
naive universality of the Wiener chaos given in \citep{bope1}, Proposition~5.4.

%th1.8 #&#
\begin{theorem}
\label{transferWignerfreepoisson}
Let $\lambda$ be a positive real number. Denote by $\hat{P}(\lambda)$
the free centered Poisson distribution with rate $\lambda$ and by
$\mathcal{S}(0,\lambda)$ the centered semicircular distribution with
variance~$\lambda$. Fix an integer $q \geq1$ and let $ \lbrace f_n
\dvtx  n \geq1  \rbrace$ be a tamed sequence of
mirror-symmetric functions in $L^2 ( \R_{+}^{q} )$ such that,
for each $n \geq1$, $\varphi (I_{q}^{\hat{N}} ( f_n )^2
 ) = \Vert f_n \Vert_{L^2 ( \R_{+}^{q} )}^{2} = \lambda$
and $\varphi (I_{q}^{\hat{N}} ( f_n )^4  ) - 2\varphi
 (I_{q}^{\hat{N}} ( f_n )^3  ) \converge2\lambda^2 -
\lambda$. Then, it holds that
\begin{enumerate}[(ii)]
\item[(i)] $I_{q}^{\hat{N}} ( f_n ) \cvlaw Z(\lambda)$ and
$I_{q}^{S} ( f_n ) \cvlaw Z(\lambda)$ if $q$ is even;
\item[(ii)] $I_{q}^{\hat{N}} ( f_n ) \cvlaw Z(\lambda)$ and
$I_{q}^{S} ( f_n ) \cvlaw\mathcal{S}(0,\lambda)$ if $q$ is odd.
\end{enumerate}
\end{theorem}

Point (i) in the above theorem implies that it suffices that the free
Poisson multiple integral (of even order) of a sequence of functions
converges towards a free Poisson limit for the Wigner integral of the
same sequence to converge to the same limit. Point (ii) provides a
counterexample to the naive universality of the Wigner chaos stated in
\citep{noudey1} since we are in a situation where the semicircular
multiple integral converges to a semicircular limit but not the free
Poisson multiple integral.

%
%re1.9 #&#
\begin{remark}
In order to avoid unnecessary heavy notations, Theorem~\ref
{fourthmomenttheoremfreepoifreepoi}, Corollary~\ref{elcorollary} and
Theorem~\ref{transferWignerfreepoisson} are stated (and proved) for
kernels in $L^2 (\R_{+}^q  ) $, for some integer $q \geq1$.
These results also hold for kernels in $L^2 ( ( \R^d )^q
 ) $, for any integer $d \geq1$, without changing anything to the proofs.
\end{remark}

%
%s1.2 #&#
\subsection{Plan}
Section~\ref{section2} provides some preliminaries on non-crossing partitions as
well as some basic definitions about the free Poisson algebra, where
our main objects of interest live. The proofs of Theorem~\ref
{fourthmomenttheoremfreepoifreepoi}, Corollary~\ref{elcorollary} and
Theorem~\ref{transferWignerfreepoisson} can be found in Section~\ref{sec3}.
Finally, Section~\ref{sec4} contains some auxiliary lemmas along with their
proofs that are used in the proofs of our main results in Section~\ref{sec3}.

%s2 #&#
\section{Preliminaries}
\label{section2}

The framework used here and the associated notation are the same as in
\citep{bope1}; see that paper for all the definitions that are not
explicitly provided here.

%s2.1 #&#
\subsection{Non-crossing partitions}
Given an integer $m \geq1$, we write $[m] = \{1,\ldots,m\}$. A \textit
{partition} of $[m]$ is a collection of non-empty and disjoint subsets
of $[m]$, called \textit{blocks}, such that their union is equal to
$[m]$. The set of all partitions of $[m]$ is denoted by $\mathcal
{P}(m)$. The cardinality of a block is called \textit{size}. We adopt
the convention of ordering the blocks of a given partition $\pi= \{
B_1,\ldots,B_r\}$ by their least element, that is: $\min B_i < \min
B_j$ if and only if $i< j$. A partition $\pi$ of $[n]$ is said to be
\textit{non-crossing} if one cannot find integers $p_1,q_1,p_2,q_2$
such that: (a) $1\leq p_1 < q_1 < p_2 < q_2\leq m$, (b) $p_1,p_2$ are
in the same block of $\pi$, (c) $q_1,q_2$ are in the same block of $\pi
$, and (d) the $p_i$'s are not in the same block of $\pi$ as the
$q_i$'s. The collection of the non-crossing partitions of $[n]$ is
denoted by $\NC(n)$, $n\geq1$. It is a well-known fact (see, e.g.,
\citep{NicSpe}, page 144) that the reversed refinement order (written
$\preceq$) induces a lattice structure on $\NC(n)$: we shall denote by
$\vee$ and $\wedge$, respectively, the associated \textit{join} and
\textit{meet} operations, where $\hat0 = \{\{1\},\ldots,\{n\}\}$ and
$\hat1=\{[n]\}$ are the corresponding minimal and maximal partitions
of the lattice. Let $m$ and $q$ be two integers such that $m,q \geq1$.
We define the partition $\pi^{*} = \{ B_1,\ldots, B_m\}\in \NC(mq)$,
where $B_1 =  \lbrace1,\ldots, q  \rbrace$, $B_2 =
\lbrace q+1,\ldots, 2q  \rbrace$ and so on until $B_m =
\lbrace(m-1)q,\ldots, mq  \rbrace$. Such a partition $\pi^{*}$ is
sometimes called a \textit{block partition}. For any integers $m,q \geq
1$, we define the four following subsets of partitions of $ [mq
 ]$:
\begin{eqnarray*}
 {\operatorname{NC}}^{0} \bigl( [mq ] , \pi^{*} \bigr) &=&
\bigl\lbrace\sigma\in{\NC}(mq) \dvtx  \sigma\wedge\pi^{*} = \hat{0} \bigr
\rbrace;
\\
 {\operatorname{NC}}^{0}_{2} \bigl( [mq ] ,
\pi^{*} \bigr) &=& \bigl\lbrace\sigma\in{\operatorname{NC}}^{0}
\bigl( [mq ] , \pi^{*} \bigr) \dvtx  \vert b \vert= 2, \forall b \in\sigma
\bigr\rbrace;
\\
 {\operatorname{NC}}^{0}_{> 2} \bigl( [mq ] ,
\pi^{*} \bigr) &=& \bigl\lbrace\sigma\in{\operatorname{NC}}^{0}
\bigl( [mq ] , \pi^{*} \bigr) \dvtx  \vert b \vert> 2, \forall b \in\sigma
\bigr\rbrace;
\\
 {\operatorname{NC}}^{0}_{\geq2} \bigl( [mq ] ,
\pi^{*} \bigr) &=& \bigl\lbrace\sigma\in{\operatorname{NC}}^{0}
\bigl( [mq ] , \pi^{*} \bigr) \dvtx  \vert b \vert\geq2, \forall b \in
\sigma \bigr\rbrace.
\end{eqnarray*}
Observe that, by definition, it holds that for any $m,q \geq1$,
${\operatorname{NC}}^{0}_{\geq2} = {\operatorname{NC}}^{0}_{2}\cup
{\operatorname{NC}}^{0}_{> 2}$ and ${\operatorname{NC}}^{0}_{2}\cap
{\operatorname{NC}}^{0}_{> 2} = \emptyset$.

Let $q,m\geq1$ be integers, and consider a function $f$ in $q$
variables. Given a partition $\sigma$ of $[mq]$, we define the function
$f_{\sigma}$, in $\vert\sigma\vert$ variables, as the mapping obtained
by identifying the variables $x_i$ and $x_j$ in the argument of the tensor
%
%e1 #&#
\begin{equation}
\label{formegeneraledefsigma} f\otimes\cdots\otimes f (x_1,\ldots,x_{mq}
) = \prod_{j=1}^{m}f (x_{(j-1)q + 1},
\ldots,x_{jq} )
\end{equation}
if and only if $i$ and $j$ are in the same block of $\sigma$.

%de2.1 #&#
\begin{definition}
Let $q\geq1$ be an integer. We say that the sequence $\{g_n \dvtx  n\geq1\}
\subset L^2(\R_{+}^q)$ is \textit{tamed} if the following conditions
hold: every $g_n$ is bounded and has bounded support and, for every
$m\geq2$ and every $\sigma\in\mathcal{P}(mq)$ such that $\sigma\wedge
\pi^* = \hat{0}$, the numerical sequence
%
%e2 #&#
\begin{equation}
\label{tamedsequencebornee} \int_{\R_{+}^{|\sigma|}} \llvert g_n \rrvert
_\sigma \,\mathrm{d}\mu^{|\sigma|},\qquad  n\geq1,
\end{equation}
is bounded, where $\pi^*\in\mathcal{P}(mq)$ is the block partition
with $m$ consecutive blocks of size $q$, and the function $\vert g_n
\vert_{\sigma}$, in $\vert\sigma\vert$ variables, is defined
according to \eqref{formegeneraledefsigma} in the case $f_n = \vert g_n
\vert$.
\end{definition}

There exists sufficient conditions in order for a sequence $\{f_n\}$ to
be tamed. It basically consists in requiring that $\{f_n\}$
concentrates asymptotically, without exploding, around a hyperdiagonal:
fix $q\geq2$, and consider a sequence $\{f_n \dvtx  n\geq1\}\subset L^2(\R
_{+}^q)$. Assume that there exist strictly positive numerical sequences
$\{M_n, z_n, \alpha_n\dvtx  n\geq1\}$ such that $\alpha_n/z_n \to0$ and
the following properties are satisfied: (a) the support of $f_n$ is
contained in the set $(-z_n,z_n)\times\cdots\times(-z_n, z_n)$
(Cartesian product of order $q$); (b) $|f_n|\leq M_n$; (c)
$f_n(x_1,\ldots,x_q) = 0$, whenever there exist $x_i, x_j$ such that
$\vert x_i - x_j\vert>\alpha_n$; (d) for every integer $m\geq q$, the
mapping $n\mapsto M_n^m z_n \alpha_n^{m-1}$ is bounded. Then, $
\lbrace f_n \dvtx  n\geq1  \rbrace$ is tamed.

%s2.2 #&#
\subsection{The free Poisson algebra}
\label{subsecfreepoissonalg}
Let $(\mathscr{A},\varphi)$ be a free tracial probability space and let
$\mathscr{A}_{+}$ denote the cone of positive operators in $\mathscr
{A}$. Denote by $\mu$ the Lebesgue measure on $\mathscr{B} ( \R
_{+} )$, where $\mathscr{B} ( \R_{+} )$ denotes the Borel
sets of $\R_{+}$ and write $\mathscr{B}_{\mu} ( \R_{+} ) =
 \lbrace B \in\mathscr{B} ( \R_{+} ) \dvtx \mu(B) <
\infty \rbrace$. The following is a brief description the free
Poisson algebra, as constructed and studied in \citep{bope1}.

For every integer $q\geq2$, the space $L^2(\R_{+}^q)$ is the
collection of all complex-valued functions on $\R_{+}^q$ that are
square-integrable with respect to $\mu^q$. Given $f\in L^2(\R_{+}^q)$,
we write $f^*(t_1,t_2,\ldots,t_q) = \overline{f(t_q,\ldots,t_2,t_1)}$,
and we call $f^*$ the \textit{adjoint} of $f$. As pointed out in
Remark~\ref{remarlmirror}, a mirror-symmetric element of $L^2(\R_{+}^q)$ is a
function $f$ that satisfies $f(t_1,\ldots,t_q) = f^*(t_1,\ldots,t_q)$,
for almost every vector $(t_1,\ldots,t_q)\in\R_{+}^q$. Observe that
mirror symmetric functions constitute a Hilbert subspace of $L^2(\R
_{+}^q)$. Let $f \in L^2 (\R_{+}^m )$ and $g \in L^2 (\R
_{+}^n )$. We define the \textit{arc} and \textit{star}
contractions of $f$ and $g$: for $1 \leq k \leq m \wedge n$, we set\vspace*{-1.5pt}
\begin{eqnarray*}
&& f \cont{k} g (t_1,\ldots, t_{m+n - 2k}) \\
&&\quad = f
\star_{k}^{k} g (t_1,\ldots, t_{m+n - 2k})
\\
&&\qquad {}\times \int_{\R_{+}^{k}}f(t_1,\ldots, t_{m-k},s_k,
\ldots, s_1)g(s_1, \ldots, s_k ,
t_{m-k+1},\ldots, t_{m+n - 2k})\mu(\mathrm{d}s_1)\cdots
\mu(\mathrm{d}s_k),
\end{eqnarray*}
and we set moreover $f \cont{0} g =f \star_{0}^{0} g = f \otimes g$.
For $k=1,\ldots,m\wedge n$, the star contraction of index $k$ (of $f$
and $g$) is defined by\vspace*{-1.5pt}
\begin{eqnarray*}
&& f \star_{k}^{k-1} g (t_1,\ldots,
t_{m+n - 2k+1}) \\
&&\quad = \int_{\R
_{+}^{k-1}}f(t_1,\ldots,
t_{m-k+1}, s_{k-1}, \ldots, s_1)
\\
&& \qquad {}\times g(s_1, \ldots, s_{k-1} , t_{m-k+1},
\ldots, t_{m+n - 2k+1})\mu (\mathrm{d}s_1)\cdots\mu(
\mathrm{d}s_{k-1}).
\end{eqnarray*}
Let $\hat{N}$ and $S$ be a \textit{free centered Poisson random
measure} and a \textit{semicircular random measure}, respectively. For
$f \in L^2 (\R_{+}^q )$, we denote by $I_{q}^{\hat{N}}(f)$
(respectively, $I_{q}^{S}(f)$) the multiple integral of $f$ with
respect to $\hat{N}$ (respectively, $S$). The space $L^2(\mathcal
{X}({\hat{N}}), \varphi) = \{I^{\hat{N}}_q(f) \dvtx  f\in L^2(\R_{+}^q),
q\geq0\}$ is a unital $\ast$-algebra, with product rule given, for
any $m,n\geq1$, $f \in L^2 (\R_{+}^m )$, $g \in L^2 (\R
_{+}^n )$, by\vspace*{-1.5pt}
%
%e3 #&#
\begin{eqnarray}
\label{e:mulp} I^{\hat{N}}_m(f)I^{\hat{N}}_n(g)
= \sum_{k=0}^{m \wedge n} I^{\hat
{N}}_{m+n-2k}
( f \cont{k} g ) + \sum_{k=1}^{m \wedge n}
I^{\hat{N}}_{m+n-2k+1} \bigl( f \star_{k}^{k-1} g
\bigr)
\end{eqnarray}
and involution $I_q^{\hat{N}}(f)^* = I_q^{\hat{N}}(f^*)$. Denote by
$L_b(\mathcal{X}(\hat{N}),\varphi)$ the collection of all objects of
the type $I_q^{\hat{N}}(f)$, where $f$ is a bounded function with
bounded support. $L_b(\mathcal{X}(\hat{N}),\varphi)$ is a subalgebra of
$L^2(\mathcal{X}(\hat{N}),\varphi)$.

The space $L^2(\mathcal{X}(S), \varphi) = \{I^{S}_q(f) \dvtx  f\in L^2(\R
_{+}^q), q\geq0\}$ is a unital $\ast$-algebra, with product rule
given, for any $m,n\geq1$, $f \in L^2 (\R_{+}^m )$, $g \in
L^2 (\R_{+}^n )$, by\vspace*{-1.5pt}
%
%e4 #&#
\begin{eqnarray}
\label{e:muls} I^S_m(f)I^S_n(g)
= \sum_{k=0}^{m \wedge n} I^S_{m+n-2k}
( f \cont {k} g )
\end{eqnarray}
and involution $I_q^{\hat{N}}(f)^* = I_q^{\hat{N}}(f^*)$.

Observe that it follows from the definition of the involution on the
algebras $L^2(\mathcal{X}(\hat{N}), \varphi)$ and $L^2(\mathcal{X}(S),
\varphi)$ that operators of the type $I_q^{\hat{N}}(f)$ or $I_q^{S}(f)$
are self-adjoint if and only if $f$ is mirror symmetric.

The following diagram formulas were proved in \citep{bope1}, Theorem~3.15, and provide an explicit combinatorial way of computing
moments of multiple integrals with respect to either a free Poisson or
a free semicircular measure: for any $f \in L^2 (\R_{+}^q )$
and any integer $q\geq1$ and $m\geq2$, it holds that
%
%e5 #&#
%e6 #&#
\begin{eqnarray}
\label{e:dp2} \varphi \bigl( I_{q}^{\hat{N}}(f)^m
\bigr) &=& \sum_{\sigma\in
{\NC}^{0}_{\geq2} ( [mq  ], \pi^{*}  ) }\int_{\R
_{+}^{\vert\sigma\vert}}f_{\sigma}
\,\mathrm{d}\mu^{\vert\sigma\vert},
\\
\label{e:ds2} \varphi \bigl( I_{q}^{S}(f)^m
\bigr) &=& \sum_{\sigma\in
{\NC}^{0}_{2} ( [mq  ], \pi^{*}  ) }\int_{\R
_{+}^{mq/2}}f_{\sigma}
\,\mathrm{d}\mu^{mq/2}.
\end{eqnarray}

%
%re2.2 #&#
\begin{remark}
These diagram formulae are free analogues of the classical diagram
formulae for multiple integrals with respect to a Gaussian or Poisson
random measure (see, e.g., \citep{PecTaq}, Theorem~7.1.3). The main
difference between the free and classical diagram formulae is the kind
of partitions of $ [mq  ]$ on which to compute the sums
appearing in the right-hand sides of the formulae. In the free case,
the sets of partitions on which to sum are subset of those of the
classical case where only non-crossing partitions are considered.
\end{remark}

%s3 #&#
\section{Proof of the main results}\label{sec3}
Let $i < m$ be two non-negative integers. For $0 \leq i \leq m-1$,
define the multisets $M_{i}^{m} =  \lbrace1,\ldots,1,0,\ldots,0
 \rbrace$ where the element $1$ has multiplicity $i$ and the
element $0$ has multiplicity $m-i-1$. Such a set is sometimes denoted
$ \lbrace(1,i),(0,m-i-1)  \rbrace$. We denote the group of
permutations of the multiset $M_{i}^{m}$ by $\mathfrak{S}_{i}^{m}$ and
its cardinality is given by the multinomial coefficient $
{m-1\choose i,m-i-1} = \frac{(m-1)!}{i!(m-i-1)!} = {m-1\choose i}$. Observe
that in the definition of the group of\vspace*{2pt} permutations of a multiset, each
permutation yields a different ordering of the elements of the
multiset, which is why the cardinality of $\mathfrak{S}_{i}^{m}$ is
${m-1\choose i}$ and not $(m-1)!$.

For any $m \geq1$, $\mathfrak{S}_{0}^{m}$ and $\mathfrak{S}_{m}^{m}$
each only have exactly one element that we denote by $\sigma_{0}$ and
$\sigma_{1}$ respectively, ($\sigma_{0}$ and $\sigma_{1}$ are in fact
the identity maps on the sets $M_{0}^{m}$ and $M_{m}^{m}$
respectively). Furthermore, for a given $\sigma\in\mathfrak
{S}_{i}^{m}$, $0 \leq i \leq m$, we define the sets
\begin{eqnarray*}
\mathfrak{A}_{m}^{\sigma} &=& \Biggl\lbrace(r_{1},
\ldots,r_{m-1}) \in (0,1,\ldots, q)^{m-1} \dvtx \forall 1 \leq p
\leq m-1,\\
&&\hphantom{\Biggl\lbrace} \sigma(p) \leq r_{p} \leq pq + \sum_{k=1}^{p-1}
\bigl(\sigma(k) - 2r_{k}\bigr) \Biggr\rbrace,
\\
\mathfrak{B}_{m}^{\sigma} &=& \Biggl\lbrace (r_{1},
\ldots,r_{m-1}) \in \mathfrak{A}_{m}^{\sigma}
\dvtx 2r_1 + \cdots+ 2r_{m-1} = mq + \sum
_{p=1}^{m-1}\sigma(p) \Biggr\rbrace,
\\
\mathfrak{D}_{m}^{\sigma} &=& \biggl\lbrace (r_{1},
\ldots,r_{m-1}) \in \mathfrak{B}_{m}^{\sigma}\cap \biggl
\lbrace0, \frac{q+1}{2}, q \biggr\rbrace \dvtx \forall1 \leq j \leq m-1,
\\
&&\hphantom{\Biggl\lbrace} r_j \in \lbrace0,q \rbrace\Leftrightarrow\sigma(j) = 0
\mbox{ and } r_j = \frac{q+1}{2} \Leftrightarrow\sigma(j) = 1
\biggr\rbrace,
\\
\mathfrak{E}_{m}^{\sigma} &=& \mathfrak{B}_{m}^{\sigma}
\setminus \mathfrak{D}_{m}^{\sigma}.
\end{eqnarray*}
In the upcoming proofs, we will drop the superscript $\hat{N}$ on free
Poisson multiple integral whenever there is only this one kind of
multiple integrals involved. Whenever a proof deals with different
sorts of multiple integrals, we will resume using the appropriate
superscripts to avoid confusion. Finally, in order to avoid more than
necessarily heavy notations, we will write
\[
_{p=1}^{m-1}\bigstar_{r_p}^{r_p - \sigma(p)}f = \bigl(
\cdots \bigl( \bigl( f \star_{r_{1}}^{r_{1} - \sigma(1)} f \bigr) \star
_{r_{2}}^{r_{2} - \sigma(2)} f \bigr) \cdots f \bigr) \star
_{r_{m-1}}^{r_{m-1} - \sigma(m-1)} f,
\]
where the $\sigma(p)$ are integers equal to either $0$ or $1$. Using
the notation introduced in Section~\ref{subsecfreepoissonalg}, we write
\[
_{p=1}^{m-1}\cont{r_p} f :=\,_{p=1}^{m-1}
\bigstar_{r_p}^{r_p}f = \bigl( \cdots \bigl( ( f
\cont{r_1} f ) \cont{r_2} f \bigr) \cdots f \bigr)
\cont{r_{m-1}}
\]
whenever all the $\sigma(p)$, $1 \leq p \leq m-1$, are zero.

%
%re3.1 #&#
\begin{remark}
The reason why the sequences of functions appearing in the main results
of the paper are required to be tamed and mirror symmetric does not
appear explicitly in the upcoming proofs. This condition ensures that
the spectral radius of free Poisson multiple integrals is bounded (see
\citep{bope1}, Theorem~3.15, for details) and that the diagram formula
\eqref{e:dp2} holds, hence guarantying the validity of the convergence
in distribution results.
\end{remark}

%s3.1 #&#
\subsection{Proof of Theorem \texorpdfstring{\protect\ref{fourthmomenttheoremfreepoifreepoi}}{1.5}}
Proving the implication (i) $\Rightarrow$ (ii)  is trivial as in the
free probability setting, convergence in distribution is equivalent to
the convergence of moments. Hence, $\varphi (I_{q}(f_n)^4  ) -
2\varphi (I_{q}(f_n)^3  ) \converge\allowbreak \varphi (Z(\lambda)^4
 ) - 2\varphi (Z(\lambda)^3  ) = 2\lambda^2 - \lambda$.
The rest of the proof will be devoted to proving the implication (ii)
$\Rightarrow$ (i). As convergence in distribution is equivalent to the
convergence of moments, we will prove that for any integer $m \geq2$,
we have $\varphi (I_{q}(f_n)^m  ) \rightarrow\varphi
(Z(\lambda)^m  )$. The proof will consist in two steps, depending
on whether $q$ is even or odd. We start with the case where $q$ is even.

\textit{Step} 1: $q$ \textit{is even}. Using Lemma~\ref{prop2}, one can write
%
%e7 #&#
\begin{eqnarray}
\label{identif2} \varphi \bigl( I_{q}(f_n)^m
\bigr) &=& \sum_{(r_{1},\ldots,r_{m-2}) \in
\mathfrak{B}_{m-1}^{\sigma_{0}}} \bigl( \cdots \bigl( (
f_n \cont {r_1} f_n ) \cont{r_2}
f_n \bigr) \cdots f \bigr) \cont{q} f_n
\nonumber
\\[-8pt]\\[-8pt]
&&{} + \sum_{i=1}^{\lfloor\vfrac{m-2}{2} \rfloor}\sum
_{\sigma\in
\mathfrak{S}_{2i}^{m-1}}\sum_{(r_{1},\ldots,r_{m-2}) \in\mathfrak
{B}_{m-1}^{\sigma}} \bigl(
_{p=1}^{m-2}\bigstar_{r_p}^{r_p - \sigma
(p)}f_n
\bigr) \cont{q}f_n.\nonumber
\end{eqnarray}
The first sum has already been addressed in \citep{nope1}, Proof of Theorem~1.4.
Under condition (ii) and by Lemma~\ref{lemma3} (note
that it is here that the special behaviour of the contraction of order
$q/2$ comes into play. See \citep{nope1}, Proof of Theorem~1.4, for
details), it holds that
\begin{eqnarray*}
\sum_{(r_{1},\ldots,r_{m-2}) \in\mathfrak{B}_{m-1}^{\sigma}} \bigl( \cdots \bigl( ( f_n
\cont{r_1} f_n ) \cont{r_2} f_n
\bigr) \cdots f_n \bigr) \cont{q} f_n \converge\varphi
\bigl(Z(\lambda)^m \bigr).
\end{eqnarray*}
Therefore, it remains to prove that
%
%e8 #&#
\begin{equation}
\label{remainstoprove} \sum_{i=1}^{\lfloor\vfrac{m-2}{2} \rfloor}\sum
_{\sigma\in\mathfrak
{S}_{2i}^{m-1}}\sum_{(r_{1},\ldots,r_{m-2}) \in\mathfrak
{B}_{m-1}^{\sigma}} \bigl(
_{p=1}^{m-2}\bigstar_{r_p}^{r_p - \sigma
(p)}f_n
\bigr) \cont{q}f_n \rightarrow0.
\end{equation}
Recalling that $\operatorname{NC}_{\geq2}^{0}$ is the disjoint union
of $\operatorname{NC}_{2}^{0}$ and $\operatorname{NC}_{>2}^{0}$ and
using the diagram formula \eqref{e:dp2}, we can write
%
%e9 #&#
\begin{equation}
\label{identif1} \varphi \bigl( I_{q}(f_n)^m
\bigr) = \sum_{\tau\in\operatorname
{NC}_{2}^{0} (  [mq  ], \pi^{*}  ) }\int_{\R_{+}^{mq
/2}}(f_n)_{\tau}
\,\mathrm{d}\mu^{mq /2} + \sum_{\tau\in\operatorname
{NC}_{>2}^{0} (  [mq  ], \pi^{*}  ) }\int
_{\R
_{+}^{\vert\tau\vert}}(f_n)_{\tau}\,\mathrm{d}
\mu^{\vert\tau\vert}.
\end{equation}
Observe that, on one hand, the diagram formula for semicircular
multiple integrals \eqref{e:ds2} states that the $m$th moment of a
semicircular multiple integral is equal to
\[
\sum_{\tau\in\operatorname{NC}_{2}^{0} (  [mq  ], \pi
^{*}  ) }\int_{\R_{+}^{mq /2}}(f_n)_{\tau}
\,\mathrm{d}\mu^{mq /2}
\]
and that on the other hand, \citep{nope1}, Proof of Theorem~1.4,
provides the following expression for it:
\[
\sum_{(r_{1},\ldots,r_{m-2}) \in\mathfrak{B}_{m-1}^{\sigma_{0}}} \bigl( \cdots \bigl( ( f_n
\cont{r_1} f_n ) \cont{r_2} f_n
\bigr) \cdots f \bigr) \cont{q} f_n.
\]
Using \eqref{identif2} and \eqref{identif1}, we get the following
identification:
\begin{eqnarray*}
 \sum_{i=1}^{\lfloor\vfrac{m-2}{2} \rfloor}\sum
_{\sigma\in\mathfrak
{S}_{2i}^{m-1}}\sum_{(r_{1},\ldots,r_{m-2}) \in\mathfrak
{B}_{m-1}^{\sigma}} \bigl(
_{p=1}^{m-2}\bigstar_{r_p}^{r_p - \sigma
(p)}f_n
\bigr) \cont{q}f_n = \sum_{\tau\in\operatorname
{NC}_{>2}^{0} (  [mq  ], \pi^{*}  ) }\int
_{\R
_{+}^{\vert\tau\vert}}(f_n)_{\tau}\,\mathrm{d}
\mu^{\vert\tau\vert}.
\end{eqnarray*}
Using the same argument as in \citep{bope1}, Proof of Theorem~4.3, it
holds that the condition $\Vert f_n \star_{r}^{r-1} f_n \Vert_{L^2
(\R_{+}^{2q - 2r +1}  )} \converge0$ for all $r \in \lbrace
1,\ldots,q  \rbrace$, implied by (ii) as stated in Lemma~\ref
{lemma3}, along with the fact that the sequence $ \lbrace f_n
\dvtx  n\geq1  \rbrace$ is tamed, is a sufficient condition in
order to have
\[
\sum_{\tau\in\operatorname{NC}_{>2} (  [mq  ], \pi^{*}
 ) }\int_{\R_{+}^{\vert\tau\vert}}(f_n)_{\tau}
\,\mathrm{d}\mu^{\vert\tau
\vert} \converge0,
\]
and hence \eqref{remainstoprove}, which concludes this step of the proof.

\textit{Step} 2: $q$ \textit{is odd}. Recall that $\mathfrak{B}_{m}^{\sigma}$ is
the disjoint union of $\mathfrak{D}_{m}^{\sigma}$ and $\mathfrak
{E}_{m}^{\sigma}$. Using Lemma~\ref{prop2}, we now have
\begin{eqnarray*}
\varphi \bigl( I_{q}(f_n)^m \bigr) &=& \sum
_{i=0}^{\lfloor\vfrac{m-2}{2}
\rfloor}\sum
_{\sigma\in\mathfrak{S}_{2i+\pi(qm)}^{m-1}}\sum_{(r_{1},\ldots,r_{m-2}) \in\mathfrak{D}_{m-1}^{\sigma}} \bigl(
_{p=1}^{m-2}\bigstar_{r_p}^{r_p - \sigma(p)}f_n
\bigr) \cont{q}f_n
\\
&&{} + \sum_{i=0}^{\lfloor\vfrac{m-2}{2} \rfloor}\sum
_{\sigma\in
\mathfrak{S}_{2i+\pi(qm)}^{m-1}}\sum_{(r_{1},\ldots,r_{m-2}) \in
\mathfrak{E}_{m-1}^{\sigma}} \bigl(
_{p=1}^{m-2}\bigstar_{r_p}^{r_p -
\sigma(p)}f_n
\bigr) \cont{q}f_n.
\end{eqnarray*}
Condition (ii), along with Lemma~\ref{lemma3}, implies that $\Vert
f_n \star_{(q+1)/2}^{(q-1)/2} f_n -f_n\Vert_{L^2 (\R_{+}^{q}
)} \converge0$ and Lemma~\ref{lemma4} ensures that, given these facts,
\begin{eqnarray*}
\sum_{i=0}^{\lfloor\vfrac{m-2}{2} \rfloor}\sum
_{\sigma\in\mathfrak
{S}_{2i+\pi(qm)}^{m-1}}\sum_{(r_{1},\ldots,r_{m-2}) \in\mathfrak
{D}_{m-1}^{\sigma}} \bigl(
_{p=1}^{m-2}\bigstar_{r_p}^{r_p - \sigma
(p)}f_n
\bigr) \cont{q}f_n \converge\varphi \bigl( Z(\lambda)^m
\bigr).
\end{eqnarray*}
It remains to show that
\begin{eqnarray*}
\sum_{i=0}^{\lfloor\vfrac{m-2}{2} \rfloor}\sum
_{\sigma\in\mathfrak
{S}_{2i+\pi(qm)}^{m-1}}\sum_{(r_{1},\ldots,r_{m-2}) \in\mathfrak
{E}_{m-1}^{\sigma}} \bigl(
_{p=1}^{m-2}\bigstar_{r_p}^{r_p - \sigma
(p)}f_n
\bigr) \cont{q}f_n \converge0.
\end{eqnarray*}
Observe that in the decomposition
\begin{eqnarray*}
\sum_{\tau\in\operatorname{NC}_{>2}^{0} (  [mq  ], \pi
^{*}  ) }\int_{\R_{+}^{\vert\tau\vert}}(f_n)_{\tau}
\,\mathrm{d}\mu^{\vert
\tau\vert}& =& \sum_{\tau\in\mathcal{C}_{>2}^{0,1} (  [mq  ], \pi^{*}
 ) }\int
_{\R_{+}^{\vert\tau\vert}}(f_n)_{\tau}\,\mathrm{d}
\mu^{\vert\tau
\vert} \\
&&{} +\sum_{\tau\in\mathcal{C}_{>2}^{0,2} (  [mq  ],
\pi^{*}  ) }\int
_{\R_{+}^{\vert\tau\vert}}(f_n)_{\tau}\,\mathrm{d}\mu
^{\vert\tau\vert},
\end{eqnarray*}
where $\mathcal{C}_{>2}^{0,1} (  [mq  ], \pi^{*}  ) =
 \lbrace\tau\in\operatorname{NC}_{>2}^{0} (  [mq
], \pi^{*}  ) \dvtx \forall b_1, b_2 \in\tau, \sharp(b_1 \cap
b_2) \in \lbrace0,\frac{q+1}{2}, q  \rbrace  \rbrace$
and $\mathcal{C}_{>2}^{0,2} (  [mq  ], \pi^{*}  ) =
\operatorname{NC}_{>2}^{0} (  [mq  ], \pi^{*}  )
\setminus\mathcal{C}_{>2}^{0,1} (  [mq  ], \pi^{*}
)$, we have
%
%e10 #&#
\begin{eqnarray}
\label{egalitecummom} &&\sum_{\tau\in\mathcal{C}_{>2}^{0,2} (  [mq  ], \pi^{*}
 ) }\int_{\R_{+}^{\vert\tau\vert}}(f_n)_{\tau}
\,\mathrm{d}\mu^{\vert\tau
\vert}\nonumber \\[-8pt]\\[-8pt]
&&\quad = \sum_{i=0}^{\lfloor\vfrac{m-2}{2} \rfloor}
\sum_{\sigma\in
\mathfrak{S}_{2i+\pi(qm)}^{m-1}}\sum_{(r_{1},\ldots,r_{m-2}) \in
\mathfrak{E}_{m-1}^{\sigma}}
\bigl( _{p=1}^{m-2}\bigstar_{r_p}^{r_p -
\sigma(p)}f_n
\bigr) \cont{q}f_n.\nonumber
\end{eqnarray}
Condition (ii) implies (through Lemma~\ref{lemma3}), that $\Vert f_n
\cont{r} f_n \Vert_{L^2 (\R_{+}^{2q - 2r}  )} \converge0$ for
all $r \in\lbrace1,\ldots,q-1 \rbrace$ and $\Vert f_n \star_{r}^{r-1}
f_n \Vert_{L^2 (\R_{+}^{2q - 2r +1}  )} \converge0$ for all
$r \in \lbrace1,\ldots,q  \rbrace\setminus \lbrace
\frac{q+1}{2} \rbrace$. As there is at least one of these
contractions appearing in each summand of the left-hand side of \eqref
{egalitecummom} (the fact that contractions appear in the left-hand
side is a direct consequence of the definition \eqref
{formegeneraledefsigma} of the quantity $(f_n)_{\tau}$ ), the argument
in \citep{bope1}, Proof of Theorem~4.3, applies once more and concludes
the proof.

%s3.2 #&#
\subsection{Proof of Corollary \texorpdfstring{\protect\ref{elcorollary}}{1.6}}
Point (i) can be proved in the same way as \citep{nope1}, Proposition~1.5, by using the contraction $f \star_{2}^{1} f$ instead of the
contraction $ f \cont{q-1} f$ in the case $q=2$.

Point (ii) can be proved by observing that: (a) if $f$ is valued in
$ \lbrace0,1  \rbrace$, then by definition of a Poisson
random measure, $I_1(f)$ has a free centered Poisson distribution with
parameter $\Vert f \Vert_{L^2(\R_+)}^{2}$, (b) if $I_1(f)$ has a free
centered Poisson distribution with parameter $\lambda> 0$, then by
Theorem~\ref{fourthmomenttheoremfreepoifreepoi} and Lemma~\ref{lemma3},
it holds that $\Vert f \star_{(q+1)/2}^{(q-1)/2} f -f\Vert_{L^2 (\R
_{+}^{q}  )} = \Vert f \star_{1}^{0} f -f \Vert_{L^2 (\R_{+}
 )} =\Vert f^2 -f \Vert_{L^2 (\R_{+}  )} = 0
\Leftrightarrow f^2 = f$, $\mu$-a.s, which concludes the proof.
\qed

%s3.3 #&#
\subsection{Proof of Theorem \texorpdfstring{\protect\ref{transferWignerfreepoisson}}{1.8}}
Point (i) is a direct consequence of Theorem~\ref
{fourthmomenttheoremfreepoifreepoi} and Theorem~1.4 in \citep{nope1}.
Point (ii) follows from the observation that when $q$ is odd, the
condition
\[
\varphi \bigl(I_{q}^{\hat{N}} ( f_n )^4
\bigr) - 2\varphi \bigl(I_{q}^{\hat{N}} ( f_n
)^3 \bigr) \converge2\lambda^2 - \lambda
\]
along with Lemma~\ref{lemma3} implies that $\Vert f_n \cont{r} f_n \Vert
_{L^2 (\R_{+}^{2q - 2r}  )} \converge0$ for all $r \in
\lbrace1,\ldots,q-1  \rbrace$. Applying Theorems 1.3 and 1.6 in
\citep{knps} ensures that $I_{q}^{S} ( f_n ) \cvlaw\mathcal
{S}(0,\lambda)$. The fact that $I_{q}^{\hat{N}} ( f_n ) \cvlaw
Z(\lambda)$ follows once again from Theorem~\ref
{fourthmomenttheoremfreepoifreepoi}.
\qed

%s4 #&#
\section{Auxiliary lemmas}\label{sec4}

%le4.1 #&#
\begin{lemma}
\label{lemma1}
Let $q \geq1$ and $m \geq2$ be integers. Let $f \in L^2(\R_{+}^q)$.
Then, it holds that
%
%e11 #&#
\begin{equation}
\label{puissancem} I_{q}(f)^m = \sum
_{i=0}^{m-1}\sum_{\sigma\in\mathfrak{S}_{i}^{m}}
\sum_{(r_{1},\ldots,r_{m-1}) \in\mathfrak{A}_{m}^{\sigma}}I_{mq+i - 2\sum
_{k=1}^{m-1}r_{k}} \bigl(
_{p=1}^{m-1}\bigstar_{r_p}^{r_p - \sigma(p)}f \bigr).
\end{equation}
\end{lemma}

\begin{pf}
The proof is done by induction on $m$. The initialization for $m = 2$
is precisely the free Poisson multiplication formula \eqref{e:mulp}.
Assume \eqref{puissancem} holds for all $p \leq m$. Then, we have
\begin{eqnarray*}
I_{q}(f)^{m+1} = \sum_{i=0}^{m-1}
\sum_{\sigma\in\mathfrak
{S}_{i}^{m}}\sum_{(r_{1},\ldots,r_{m-1}) \in\mathfrak{A}_{m}^{\sigma
}}I_{mq+i - 2\sum_{k=1}^{m-1}r_{k}}
\bigl( _{p=1}^{m-1}\bigstar _{r_p}^{r_p - \sigma(p)}f
\bigr)I_{q}(f).
\end{eqnarray*}
We use the multiplication formula \eqref{e:mulp} once again to obtain
\begin{eqnarray*}
&& I_{q}(f)^{m+1}
\\
&& \quad = \sum_{i=0}^{m-1}\sum
_{\sigma\in\mathfrak{S}_{i}^{m}}\sum_{(r_{1},\ldots,r_{m-1}) \in\mathfrak{A}_{m}^{\sigma}}\sum
_{r_m =
0}^{q \wedge [mq + i - 2\sum_{k=1}^{m-1}r_k  ] }I_{(m+1)q+i -
2\sum_{k=1}^{m}r_{k}} \bigl( \bigl(
_{p=1}^{m-1}\bigstar_{r_p}^{r_p -
\sigma(p)}f \bigr)
\cont{r_p} f \bigr)
\\
&&\qquad {} + \sum_{i=0}^{m-1}\sum
_{\sigma\in\mathfrak{S}_{i}^{m}}\sum_{(r_{1},\ldots,r_{m-1}) \in\mathfrak{A}_{m}^{\sigma}}\sum
_{r_m =
1}^{q \wedge [mq + i - 2\sum_{k=1}^{m-1}r_k  ] }
\\
&&\hspace*{151pt} I_{(m+1)q+(i+1) - 2\sum_{k=1}^{m}r_{k}} \bigl( \bigl( _{p=1}^{m-1}
\bigstar_{r_p}^{r_p - \sigma(p)}f \bigr) \star_{r_p}^{r_p
-1}
f \bigr)
\end{eqnarray*}
We now write the first summand of the first sum and the last summand of
the second sum separately, as they have to be treated differently from
the other terms. This yields
\begin{eqnarray*}
&& I_{q}(f)^{m+1}
\\
&&\quad  = \sum_{(r_{1},\ldots,r_{m-1}) \in\mathfrak{A}_{m}^{\sigma_{0}}}\sum_{r_m = 0}^{q \wedge [mq - 2\sum_{k=1}^{m-1}r_k  ] }I_{(m+1)q
- 2\sum_{k=1}^{m}r_{k}}
\bigl( _{p=1}^{m}\cont{r_p}f \bigr)
\\
&&\qquad {} + \sum_{i=1}^{m-1}\sum
_{\sigma\in\mathfrak{S}_{i}^{m}}\sum_{(r_{1},\ldots,r_{m-1}) \in\mathfrak{A}_{m}^{\sigma}}\sum
_{r_m =
0}^{q \wedge [mq +\sum_{k=1}^{m-1}(\sigma(k) - 2r_k)  ] }
\\
&&\hspace*{151pt} I_{(m+1)q+i - 2\sum_{k=1}^{m}r_{k}} \bigl( \bigl( _{p=1}^{m-1}
\bigstar_{r_p}^{r_p - \sigma(p)}f \bigr) \cont{r_p} f \bigr)
\\
&&\qquad {} + \sum_{i=0}^{m-2}\sum
_{\sigma\in\mathfrak{S}_{i}^{m}}\sum_{(r_{1},\ldots,r_{m-1}) \in\mathfrak{A}_{m}^{\sigma}}\sum
_{r_m =
1}^{q \wedge [mq +\sum_{k=1}^{m-1}(\sigma(k) -2r_k)  ] }
\\
&& \hspace*{151pt}I_{(m+1)q+(i+1) - 2\sum_{k=1}^{m}r_{k}} \bigl( \bigl( _{p=1}^{m-1}
\bigstar_{r_p}^{r_p - \sigma(p)}f \bigr) \star_{r_p}^{r_p
-1}
f \bigr)
\\
&&\qquad {} + \sum_{(r_{1},\ldots,r_{m-1}) \in\mathfrak{A}_{m}^{\sigma_{1}}}\sum_{r_m = 1}^{q \wedge [mq -2\sum_{k=1}^{m-1}2r_k + (m-1)  ]
}I_{(m+1)q+m - 2\sum_{k=1}^{m}r_{k}}
\bigl( _{p=1}^{m}\bigstar _{r_p}^{r_p - 1}f
\bigr).
\end{eqnarray*}
Also remark that the $i$ appearing in the upper limit of the sum on
$r_m$ has been replaced by $\sum_{k=1}^{m-1}\sigma(k)$ (according to
the definition of the sets $\mathfrak{A}_{m}^{\sigma}$). The next step
consists of shifting the index up in the third sum above. We obtain
\begin{eqnarray*}
&& I_{q}(f)^{m+1}
\\
&&\quad  = \sum_{(r_{1},\ldots,r_{m-1}) \in\mathfrak{A}_{m}^{\sigma_{0}}}\sum_{r_m = 0}^{q \wedge [mq - 2\sum_{k=1}^{m-1}r_k  ] }I_{(m+1)q
- 2\sum_{k=1}^{m}r_{k}}
\bigl( _{p=1}^{m}\cont{r_p}f \bigr)
\\
&&\qquad {} + \sum_{i=1}^{m-1}\sum
_{\sigma\in\mathfrak{S}_{i}^{m}}\sum_{(r_{1},\ldots,r_{m-1}) \in\mathfrak{A}_{m}^{\sigma}}\sum
_{r_m =
0}^{q \wedge [mq +\sum_{k=1}^{m-1}(\sigma(k) - 2r_k)  ] }
\\
&&\hspace*{151pt} I_{(m+1)q+i - 2\sum_{k=1}^{m}r_{k}} \bigl( \bigl( _{p=1}^{m-1}
\bigstar_{r_p}^{r_p - \sigma(p)}f \bigr) \cont{r_p} f \bigr)
\\
&&\qquad {} + \sum_{i=1}^{m-1}\sum
_{\sigma\in\mathfrak{S}_{i+1}^{m}}\sum_{(r_{1},\ldots,r_{m-1}) \in\mathfrak{A}_{m}^{\sigma}}\sum
_{r_m =
1}^{q \wedge [mq +\sum_{k=1}^{m-1}(\sigma(k) -2r_k)  ] }
\\
&&\hspace*{157pt} I_{(m+1)q+i - 2\sum_{k=1}^{m}r_{k}} \bigl( \bigl( _{p=1}^{m-1}
\bigstar_{r_p}^{r_p - \sigma(p)}f \bigr) \star_{r_p}^{r_p
-1}
f \bigr)
\\
&&\qquad {} + \sum_{(r_{1},\ldots,r_{m-1}) \in\mathfrak{A}_{m}^{\sigma_{1}}}\sum_{r_m = 1}^{q \wedge [mq -2\sum_{k=1}^{m-1}2r_k + (m-1)  ]
}I_{(m+1)q+m - 2\sum_{k=1}^{m}r_{k}}
\bigl( _{p=1}^{m}\bigstar _{r_p}^{r_p - 1}f
\bigr).
\end{eqnarray*}
Recalling the definitions of $\mathfrak{S}_{i}^{m}$ and $\mathfrak
{A}_{m}^{\sigma}$, one can combine the two middle sums into a single
one in the following way:
\begin{eqnarray*}
&& I_{q}(f)^{m+1}
\\
&&\quad  = \sum_{(r_{1},\ldots,r_{m-1}) \in\mathfrak{A}_{m}^{\sigma_{0}}}\sum_{r_m = 0}^{q \wedge [mq - 2\sum_{k=1}^{m-1}r_k  ] }I_{(m+1)q
- 2\sum_{k=1}^{m}r_{k}}
\bigl( _{p=1}^{m}\cont{r_p}f \bigr)
\\
&&\qquad {} + \sum_{i=1}^{m-1}\sum
_{\sigma\in\mathfrak{S}_{i}^{m+1}}\sum_{(r_{1},\ldots,r_{m}) \in\mathfrak{A}_{m+1}^{\sigma}}
I_{(m+1)q+i -
2\sum_{k=1}^{m}r_{k}} \bigl( _{p=1}^{m}\bigstar_{r_p}^{r_p - \sigma(p)}f
\bigr)
\\
&&\qquad {} + \sum_{(r_{1},\ldots,r_{m-1}) \in\mathfrak{A}_{m}^{\sigma_{1}}}\sum_{r_m = 1}^{q \wedge [mq -2\sum_{k=1}^{m-1}2r_k + (m-1)  ]
}I_{(m+1)q+m - 2\sum_{k=1}^{m}r_{k}}
\bigl( _{p=1}^{m}\bigstar _{r_p}^{r_p - 1}f
\bigr).
\end{eqnarray*}
It remains to combine the three final sums to conclude the proof.
\end{pf}

%le4.2 #&#
\begin{lemma}
\label{prop2}
Let $q \geq1$ and $m \geq2$ be integers. Let $f \in L^2(\R_{+}^q)$.
Then, it holds that
%
%e12 #&#
\begin{equation}
\label{traceofmultipleintegraltothepowerm} \varphi \bigl( I_{q}(f)^m \bigr) = \sum
_{i=0}^{\lfloor\vfrac{m-2}{2}
\rfloor}\sum
_{\sigma\in\mathfrak{S}_{2i+\pi(qm)}^{m-1}}\sum_{(r_{1},\ldots,r_{m-2}) \in\mathfrak{B}_{m-1}^{\sigma}} \bigl(
_{p=1}^{m-2}\bigstar_{r_p}^{r_p - \sigma(p)}f \bigr)
\cont{q}f,
\end{equation}
where $\pi$ is the parity function defined on $\mathbb{N}$ by $\pi(x) =
0$ if $x$ is even and $1$ otherwise.
\end{lemma}

\begin{pf}
When using Lemma~\ref{lemma1} to evaluate $\varphi (I_{q}(f)^m
)$, one has to determine when $mq - 2\sum_{k=1}^{m-1}r_{k} +\sum_{k=1}^{m-1}\sigma(k)$ is zero. In order to do so, we will study the
quantity $r_1 + \cdots+ r_{m-1}$ and determine the sufficient and
necessary conditions for it to be equal to $mq + \sum_{k=1}^{m-1}\sigma
(k)$. Recall that $(r_{1},\ldots,r_{m-1}) \in\mathfrak{A}_{m}^{\sigma
}$ and set $\zeta_{p} = (p+1)q + \sum_{k=1}^{p-1}\sigma(k)$. We will
proceed by induction to prove that, for all $p \leq m-1 $,
%
%e13 #&#
\begin{equation}
\label{inductionhypo} 2\max_{(r_{1},\ldots,r_{p}) \in\mathfrak{A}_{p+1}^{\sigma}} (r_1 + \cdots+
r_p ) = \lleft\{ %
\begin{array} {l@{\qquad}l} \zeta_{p}&
\mbox{if } \zeta_{p} \mbox{ is even},
\\
\zeta_{p} - 1& \mbox{if } \zeta_{p} \mbox{ is odd}.
\end{array} %
\rright.
\end{equation}
For $p = 1$, it is obvious that $2\max_{\sigma(1) \leq r_{1} \leq
q} r_1 = 2q$. Fix $p \leq m-2$ and assume \eqref{inductionhypo} is
verified up to rank $p$. Using the induction hypothesis, it is easy to
verify that, for $q \geq2$,
\[
2\max_{(r_{1},\ldots,r_{p}) \in\mathfrak{A}_{p+1}^{\sigma}} (r_1 + \cdots+ r_p )
\geq pq + \sum_{k=1}^{p}\sigma(k).
\]
We know that, on $r_{p+1}$, we have the restriction
\[
\sigma(p+1) \leq r_{p+1} \leq q \wedge \Biggl((p+1)q - 2\sum
_{k=1}^{p}r_k + \sum
_{k=1}^{p}\sigma(k) \Biggr).
\]
Hence, if $q \geq(p+1)q - 2\sum_{k=1}^{p}r_k + \sum_{k=1}^{p}\sigma(k)
\Leftrightarrow r_1 + \cdots+ r_p \geq\frac{pq + \sum_{k=1}^{p}\sigma
(k)}{2}$, then, if $pq + \sum_{k=1}^{p}\sigma(k) = \zeta_{p+1}$ is even,
\begin{eqnarray*}
&& 2\max_{(r_{1},\ldots,r_{p+1}) \in\mathfrak{A}_{p+2}^{\sigma
}} (r_1 + \cdots+ r_{p+1}
) \\
&&\quad = (p+1)q -2 \min_{r_{1} +
\cdots+ r_{p}}\sum_{k=1}^{p}r_k
+ \sum_{k=1}^{p}\sigma(k)
\\
&&\quad  = 2(p+1)q -pq - \sum_{k=1}^{p}\sigma(k)
+ 2\sum_{k=1}^{p}\sigma(k) \\
&&\quad = (p+2)q + \sum
_{k=1}^{p}\sigma(k).
\end{eqnarray*}
If $pq + \sum_{k=1}^{p}\sigma(k) = \zeta_{p+1}$ is odd, then
\begin{eqnarray*}
 2\max_{(r_{1},\ldots,r_{p+1}) \in\mathfrak{A}_{p+2}^{\sigma
}} (r_1 + \cdots+ r_{p+1}
) &=& (p+1)q -2 \min_{r_{1} +
\cdots+ r_{p}}\sum_{k=1}^{p}r_k
+ \sum_{k=1}^{p}\sigma(k)
\\
& =& 2(p+1)q -pq - \sum_{k=1}^{p}\sigma(k)
-1 + 2\sum_{k=1}^{p}\sigma (k)\\
& =& (p+2)q +
\sum_{k=1}^{p}\sigma(k) -1.
\end{eqnarray*}
It now remains to consider the case where $q \leq(p+1)q - 2\sum_{k=1}^{p}r_k + \sum_{k=1}^{p}\sigma(k) \Leftrightarrow r_1 + \cdots+
r_p \leq\frac{pq + \sum_{k=1}^{p}\sigma(k)}{2}$.

If $pq + \sum_{k=1}^{p}\sigma(k) = \zeta_{p+1}$ is even,
\begin{eqnarray*}
2\max_{(r_{1},\ldots,r_{p+1}) \in\mathfrak{A}_{p+2}^{\sigma}} (r_1 + \cdots+ r_{p+1} ) =
2q + 2 \max_{r_{1} + \cdots+
r_{p}}\sum_{k=1}^{p}r_k
= (p+2)q + \sum_{k=1}^{p}\sigma(k).
\end{eqnarray*}
If $pq + \sum_{k=1}^{p}\sigma(k) = \zeta_{p+1}$ is odd, then
\begin{eqnarray*}
2\max_{(r_{1},\ldots,r_{p+1}) \in\mathfrak{A}_{p+2}^{\sigma}} (r_1 + \cdots+ r_{p+1} ) =
2q + 2 \max_{r_{1} + \cdots+
r_{p}}\sum_{k=1}^{p}r_k
= (p+2)q + \sum_{k=1}^{p}\sigma(k) -1.
\end{eqnarray*}
This completes the induction.

Coming back to finding the necessary and sufficient conditions in order
to have
$2\sum_{k=1}^{m-1}r_{k} = mq + \sum_{k=1}^{m-1}\sigma(k)$, the above
result for $p = m-1$ shows that it is necessary that $\zeta_m$ is even
and that $\sigma(m-1) = 0$ for this equality to hold.

Note that if $q$ is even, then it suffices that $\sum_{k=1}^{m-2}\sigma
(k)$ be even as well for $\zeta_m$ to be even. In this case, as $\sigma
(m-1)$ has to be zero, it implies that $\sum_{k=1}^{m-1}\sigma(k)$ is
even as well. This only happens on the groups of permutations with an
even index such as $\mathfrak{S}_{2i}$.
Finally, because $2\max_{(r_{1},\ldots,r_{p}) \in\mathfrak
{A}_{p+1}^{\sigma}} (r_1 + \cdots+
r_{m-2} ) = \zeta_{m-1}$, it forces $r_{m-1}$ to always be equal
to $q$. As $\sigma(m-1)$ is always $0$ and $r_{m-1}$ is always $q$,
there is in fact no sum on $r_{m-1}$ anymore and the groups of
permutations that have to appear in \eqref
{traceofmultipleintegraltothepowerm} need only be the ones on sets of
size $m-2$. Combining these conditions yields the desired result.

It remains to examine the case where $q$ is odd. In this case, if $m$
is even, then $\sum_{k=1}^{m-2}\sigma(k)$ has to be even as well in
order for $\zeta_m$ to be even and the same arguments as in the
previous case apply. The only (slightly) different case is whenever $q$
and $m$ are odd. In this case, $\sum_{k=1}^{m-2}\sigma(k)$ has to be
odd as well in order for $\zeta_m$ to be even, and one has to consider
the groups of permutations with an odd index such as $\mathfrak
{S}_{2i+1}$ instead of the groups $\mathfrak{S}_{2i}$ for the announced
result to follow. This proves that the parity of the groups of
permutations to consider has to be the same as the parity of the
product $qm$.
\end{pf}

%%%%%%%%%%%%%%%%%%%%%%%%%%%%%%%%%%%%%%%%%%%%%%%%%%%%%%%%%%%%%%%
%
%le4.3 #&#
\begin{lemma}
\label{lemma3}
Let $q \geq1$ be an integer, and consider a sequence of functions
$ \lbrace f_n \dvtx  n \geq1  \rbrace\subset L^2 ( \R
_+^q ) $ such that $\Vert f_n \Vert_{L^2 (\R_{+}^{q}
)}^{2} = \lambda> 0 $ for every $n \geq1$. Then,
\begin{eqnarray*}
\varphi \bigl(I_{q}(f_n)^4 \bigr) - 2\varphi
\bigl(I_{q}(f_n)^3 \bigr) \converge2
\lambda^2 - \lambda
\end{eqnarray*}
\begin{enumerate}[(ii)]
\item[(i)] if and only if $\Vert f_n\cont{q/2} f_n -f_n\Vert_{L^2
(\R_{+}^{q}  )} \converge0$, $\Vert f_n \cont{r} f_n \Vert
_{L^2 (\R_{+}^{2q - 2r}  )} \converge0$ for all $r \in
\lbrace1,\ldots,q-1  \rbrace\setminus \lbrace\frac
{q}{2} \rbrace$, and $\Vert f_n \star_{r}^{r-1} f_n \Vert_{L^2
(\R_{+}^{2q - 2r +1}  )} \converge0$ for all $r \in \lbrace
1,\ldots,q  \rbrace$ if $q$ is even;
\item[(ii)] if and only if $\Vert f_n \star_{(q+1)/2}^{(q-1)/2} f_n
-f_n\Vert_{L^2 (\R_{+}^{q}  )} \converge0$, $\Vert f_n \cont
{r} f_n \Vert_{L^2 (\R_{+}^{2q - 2r}  )} \converge0$ for all
$r \in \lbrace1,\ldots,q-1  \rbrace$, and $\Vert f_n \star
_{r}^{r-1} f_n \Vert_{L^2 (\R_{+}^{2q - 2r +1}  )} \converge
0$ for all $r \in \lbrace1,\ldots,q  \rbrace\setminus
\lbrace\frac{q+1}{2} \rbrace$ if $q$ is odd;
\end{enumerate}
\end{lemma}

\begin{pf}
Compared to the proof of \citep{nope1}, Lemma~5.1, only the case where
$q$ is odd differs slightly. In that case, the product formula \eqref
{e:mulp} and orthogonality in $L^2 (\mathscr{A},\varphi ) $
of multiple integrals of different orders yield
\begin{eqnarray*}
\varphi \bigl( I_{q}(f_n)^2 -
I_{q}(f_n) \bigr) &=& 2\lambda^2 + \bigl\Vert
f_n \star_{(q+1)/2}^{(q-1)/2} f_n
-f_n\bigr\Vert_{L^2 (\R_{+}^{q}
)}^{2} + \sum
_{r=1}^{q-1} \Vert f_n \cont{r}
f_n \Vert_{L^2 (\R
_{+}^{2q - 2r}  )}^2
\\
&&{} + \mathop{\sum_{ 1 \leq r \leq q}}_{r \neq(q+1)/2}\bigl\Vert
f_n \star_{r}^{r-1} f_n
\bigr\Vert_{L^2 (\R_{+}^{2q - 2r +1}  )}^2.
\end{eqnarray*}
The conclusion is obtained as in the proof of \citep{nope1}, Lemma~5.1.
\end{pf}

%%%%%%%%%%%%%%%%%%%%%%%%%%%%%%%%%%%%%%%%%%%%%%%%%%%%%%%%%%%%%%%
%
%le4.4 #&#
\begin{lemma}
\label{lemma4}
Let $q\geq1$ be an odd integer and let $m \geq2$ be an integer. Let
$ \lbrace f_n \dvtx  n \geq1  \rbrace\subset L^2 ( \R
_{+}^{q} )$ be a sequence of tamed mirror symmetric functions such
that $\Vert f_n \Vert_{L^2 ( \R_{+}^{q} )}^{2} = \lambda>0$
for every $n$. Then, if $\Vert f _n\star_{(q+1)/2}^{(q-1)/2} f_n
-f_n\Vert_{L^2 (\R_{+}^{q}  )} \converge0$, then
\begin{eqnarray*}
\sum_{i=0}^{\lfloor\vfrac{m-2}{2} \rfloor}\sum
_{\sigma\in\mathfrak
{S}_{2i+\pi(qm)}^{m-1}}\sum_{(r_{1},\ldots,r_{m-2}) \in\mathfrak
{D}_{m-1}^{\sigma}} \bigl(
_{p=1}^{m-2}\bigstar_{r_p}^{r_p - \sigma
(p)}f_n
\bigr) \cont{q}f_n \converge\varphi \bigl( Z(\lambda)^m
\bigr) = \sum_{j=1}^{m}\lambda^j
R_{m,j},
\end{eqnarray*}
where $R_{m,j}$ is the number of non-crossing
partitions of $ [ m ] $ with exactly $j$ blocks and with no
singletons. Notice that, when $m$ is even, one has that $R_{m,j} = 0$
for every $j > m/2$ and when $m$ is odd, then $R_{m,j} = 0$ for every
$j > (m - 1)/2$. The numbers $R_{m,j}$ are related to the so-called
Riordan numbers $ \lbrace R_m \dvtx  m\geq1 \rbrace$ (for a
detailed combinatorial analysis of these numbers, see \citep{bern1}) by
$R_m = \sum_{j=1}^{m} R_{m,j}$ for all $m \geq1$).
\end{lemma}

\begin{pf}
The same arguments as in the proof of \citep{nope1}, Lemma~5.2, can be
used by replacing the case $q=2$ in the last part (where the argument
of two polynomials coinciding on a countable set being necessarily
equal is used) by the case where $q=1$ with a sequence $f_n = f = \sum_{i=1}^{p}\mathds{1}_{A_i}$, where $ \lbrace A_i\dvtx  i = 1,\ldots
,p  \rbrace$ are disjoint Borel sets with measure 1.
\end{pf}

%\begin{appendix}
%\section{}
%\end{appendix}

% zodis "Acknowledgments" paliekamas pagal autoriu
%\section*{Acknowledgements}

%\begin{supplement}%[id=suppA]
%\sname{Supplement A}
%\stitle{}
%\slink[doi]{10.3150/00-BEJXXXXSUPP} %[doi,text={...}] - jei reikia
%suskaldyti doi
%\sdatatype{.pdf}
%\sfilename{BEJ000\_supp.pdf}
%\sdescription{}
%\end{supplement}

% imsref loaded by jurgita.kaciuliene, 2014-07-01 09:00:05

\printhistory
\end{document}